\documentstyle[epsf]{article}

\newcommand{\bigzerou}{%
\smash{\lower1.7ex\hbox{\bg 0}}}
\newtheorem{theorem}{Theorem}

\newtheorem{defi}{Definition}
\newtheorem{cor}{Corollary}
\newtheorem{conj}{Conjecture}

\newcommand{\ba}{\begin{eqnarray}}
\newcommand{\ea}{\end{eqnarray}}
\newcommand{\no}{\nonumber}

\newcommand{\mapright}[1]{%
\smash{\mathop{%
\hbox to 1.0cm{\rightarrowfill}}\limits^{#1}}}
\newcommand{\mapleft}[1]{%
\smash{\mathop{%
\hbox to 1.3cm{\leftarrowfill}}\limits^{#1}}}

\begin{document}
\title{
\begin{flushright}
  \begin{minipage}[b]{5em}
    \normalsize
    ${}$      \\
  \end{minipage}
\end{flushright}
{\bf Virtual Structure Constants as Intersection Numbers of Moduli
 Space of Polynomial Maps with Two Marked Points}}
\author{Masao Jinzenji\\
\\
\it Division of Mathematics, Graduate School of Science \\
\it Hokkaido University \\
\it  Kita-ku, Sapporo, 060-0810, Japan\\
{\it e-mail address: jin@math.sci.hokudai.ac.jp}}
\maketitle
\begin{abstract}
In this paper, we derive the virtual structure constants used in the mirror
 computation of the degree $k$ hypersurface in $CP^{N-1}$,
 by using a localization computation applied to moduli space of
 polynomial maps from $CP^{1}$ to $CP^{N-1}$ with two marked points. 
We also apply this technique to the non-nef local geometry ${\cal
 O}(1)\oplus {\cal O}(-3)\rightarrow CP^{1}$ and realize the
mirror computation without using Birkhoff factorization.
\end{abstract}
\section{Introduction}
Analysis of mirror symmetry from the point of view of the Gauged Linear
 Sigma Model is very important both in the mathematical 
and physical study of mirror symmetry \cite{witten},\cite{mp}. As was
 suggested by Morrison and Plesser in \cite{mp}, the
Gauged Linear Sigma Model is directly connected to the B-model in mirror
 symmetry. Indeed, they constructed the moduli space of instantons 
of the Gauged Linear Sigma Models corresponding to Calabi-Yau 3-folds in
 $P(1,1,1,1,1)(=CP^{4})$ and (the blow-up of)$P(1,1,2,2,2)$, 
computed intersection numbers on these moduli spaces and showed that
 their generating functions coincide with  the 
B-model Yukawa couplings used in the mirror computation. Further analysis
 of mirror symmetry of Calabi-Yau 3-folds in this direction was
also pursued by Batyrev and Materov \cite{baty}. 

In this paper, we try to generalize this kind of analysis to the mirror
 computations 
of non-nef geometries. Mirror symmetry of non-nef geometries has been
 studied by various 
authors,\cite{blly},\cite{givc},\cite{iri},\cite{gene},\cite{fj1},\cite{fj2}. 
Since the mirror computation in this case is rather complicated, it is hard
 to define objects that are 
directly connected to the B-model, or to the Gauged Linear Sigma Model. In this
 paper, we compute the "virtual structure constants" 
$\tilde{L}_{n}^{N,k,d}$ used in mirror computation of a non-nef degree
 $k$ hypersurface in $CP^{N-1}$ $(k>N)$ \cite{gene}. Originally, 
the virtual structure constants are defined by recursive formulas that
 represent $\tilde{L}_{n}^{N,k,d}$ in terms of weighted 
homogeneous polynomials in $\tilde{L}_{m}^{N+1,k,d'},\;\; (d'\leq
 d)$\cite{fano}. Later, we showed that $\tilde{L}_{n}^{N,k,d}$ can be 
directly computed from the virtual Gauss-Manin system associated with
 Givental's ODE: 
\begin{equation}
\biggl((\frac{d}{dx})^{N-1}-k\cdot e^{x}\cdot
 (k\frac{d}{dx}+k-1)(k\frac{d}{dx}+k-2)\cdots(k\frac{d}{dx}+1)\biggr)w(x)=0,
\label{eq1}
\end{equation}
\cite{gm}. In \cite{iri}, Iritani showed that the virtual structure
 constants can be obtained after Birkhoff factorization of 
connection matrices of the Gauss-Manin system associated with Givental's
 ODE with "$h$" parameter:
\begin{equation}
\biggl(h^{N-1}(\frac{d}{dx})^{N-1}-h^{k-1}k\cdot e^{x}\cdot
 (k\frac{d}{dx}+(k-1))(k\frac{d}{dx}+(k-2))\cdots(k\frac{d}{dx}+1)\biggr)w(x,h)=0,
\label{eq2}
\end{equation}
generalizing the method invented by Guest et al
 \cite{guest},\cite{sakai}.

In this paper, we start from the recursive formula in \cite{gm} that
 determines $\tilde{L}_{n}^{N,k,d}$. Next, we propose a
"conjectural" residue integral representation of
 $\tilde{L}_{n}^{N,k,d}$, which is speculated by solving the recursive formula
for low degrees. Then we show that this representation can be
 interpreted as the result obtained by applying the localization 
computation to the moduli space of polynomial maps of degree $d$ with two
 marked points. By the word "polynomial map", we mean a birational map $p$
of degree $d$ from $CP^{1}$ to $CP^{N-1}$ given by,
\begin{equation}
p(s:t)=( \sum_{j=0}^{d}a_{j}^{1}s^{j}t^{d-j} :
 \sum_{j=0}^{d}a_{j}^{2}s^{j}t^{d-j}:\cdots : \sum_{j=0}^{d}a_{j}^{N}s^{j}t^{d-j} ).
\label{poly}
\end{equation}
The moduli space of polynomial maps of degree $d$ can be identified with
 $CP^{N(d+1)-1}$, which is the moduli space of instantons of the
Gauged Linear Sigma Model. The two marked points are fixed to $(0:1),\;
 (1:0)\in CP^{1}$. Then we consider polynomial maps of degree $d$ 
such that the image of two marked points are well-defined in $CP^{N-1}$ and
 divide the corresponding moduli space by ${\bf C}^{\times}$, the automorphism 
group of $CP^{1}$ fixing the two marked points. After resolving the
 singularities of the resulting space, 
we obtain the moduli space of polynomial maps with two marked points
 mentioned above. With this set up, we can derive 
the residue integral representation by applying the localization
 computation to this space. 

Our geometrical derivation using localization has the following by-product. We can
 also apply this technique to the local geometries 
$\oplus_{j=1}^{m}{\cal O}(k_{j})\rightarrow CP^{N-1}$. In this paper,
 we consider two examples of local geometry, 
${\cal O}(-1)\oplus{\cal O}(-1)\rightarrow CP^{1}$ and  ${\cal
 O}(1)\oplus{\cal O}(-3)\rightarrow CP^{1}$. 
By applying the localization computation to the moduli space of polynomial maps
 with two marked points, we can compute the
"virtual structure constants" for these models. In the latter case, it
 was hard for us to define virtual structure 
constants because we don't have an appropriate Picard-Fuchs equation like
 (\ref{eq1}). These virtual structure constants 
give us the mirror map and B-model-like two point functions, which is
 expected from analogy with the behavior of the virtual
structure constants of the quintic 3-fold in $CP^{4}$ \cite{cj}. As a
 result, we can perform the mirror computation of 
${\cal O}(1)\oplus{\cal O}(-3)\rightarrow CP^{1}$ without using
 Birkhoff factorization that had been inevitable 
in our previous analysis \cite{fj1}. 

This paper is organized as follows. In Section 2, we briefly introduce
 the virtual structure constants $\tilde{L}_{n}^{N,k,d}$ and 
its residue integral representation. In Section 3, we define the moduli
 space of polynomial maps of degree $d$ with two marked points 
, introduce a torus action on this space and determine the fixed point sets
 under the torus action. Then we derive the residue integral 
representation introduced in Section 2 by using a localization
 computation. In Section 4, we apply the method of Section 3 to 
${\cal O}(-1)\oplus{\cal O}(-1)\rightarrow CP^{1}$ and  ${\cal
 O}(1)\oplus{\cal O}(-3)\rightarrow CP^{1}$. 
Section 5 gives concluding remarks.
\\
\\
{\bf Acknowledgment} The author would like to thank Prof. D.Matsushita
 and Dr. B.Forbes for valuable discussions. 
He would also like to thank Miruko Jinzenji for encouraging him to keep
 attention to the line of thought used in this paper.   
\section{Virtual Structure Constants}
\subsection{Virtual Structure Constants and Givental's ODE}
In this subsection, we introduce the virtual structure constants 
$\tilde{L}^{N,k,d}_{m}$ which are non-zero only if $0\leq m\leq
 N-1+(k-N)d$.
The original definition of $\tilde{L}^{N,k,d}_{m}$ in \cite{gm} 
is given by the initial condition:
\begin{equation}
\sum_{m=0}^{k-1}\tilde{L}_{m}^{N,k,1}w^{m}=
k\cdot\prod_{j=1}^{k-1}(jw+(k-j)),\;\;(N-k\geq 2),
\label{vi}
\end{equation}
and the recursive formulas that describe $\tilde{L}^{N,k,d}_{m}$ as 
a weighted homogeneous polynomial in
 $\tilde{L}^{N+1,k,d'}_{n}\;\;(d'\leq d)$
of degree $d$. See \cite{fano} 
for the explicit form of the recursive formulas. In \cite{gm}, we
 showed that the virtual structure constants are directly connected with
 Givental's
ODE:     
\begin{equation}
\biggl((\frac{d}{dx})^{N-1}-k\cdot e^{x}\cdot
 (k\frac{d}{dx}+k-1)(k\frac{d}{dx}+k-2)\cdots(k\frac{d}{dx}+1)\biggr)w(x)=0,
\label{fun}
\end{equation}
for arbitrary $N$ and $k$ via the virtual Gauss-Manin system defined
as follows: 
\begin{defi}
We call the following rank 1 ODE for the vector valued function 
$\tilde{\psi}_{m}(x)$
\begin{eqnarray}
\frac{d\tilde{\psi}_{N-2-m}(x)}{dx}&=&\tilde{\psi}_{N-1-m}(x)+
\sum_{d=1}^{\infty}
\exp(dx)\cdot \tilde{L}_{m}^{N,k,d}\cdot\tilde{\psi}_{N-1-m-(N-k)d}(x).
\label{gm1}
\end{eqnarray}
the virtual Gauss-Manin system associated with the quantum K\"ahler 
sub-ring of $M_{N}^{k}$, where $m$ runs through $0\leq m \leq N-2$ 
if $N-k \geq 1$, $0\leq m \leq N-1$ if $N-k=0$, and $m\in {\bf Z}$ if 
$N-k<0$.  
\end{defi}
Here, we restate the main result in \cite{gm}.
\begin{theorem}
We can derive the following relation from the virtual Gauss-Manin
 system 
(\ref{gm1}).
\begin{equation}
\tilde{\psi}_{N-1}(x)=\biggl((\frac{d}{dx})^{N-1}-k\cdot e^{x}\cdot
(k\frac{d}{dx}+k-1)\cdots
(k\frac{d}{dx}+2)\cdot(k\frac{d}{dx}+1)\biggr)\bigl(\frac{d}{dx}\bigr)^{\beta}
\tilde{\psi}_{-\beta}(x)
\label{hyper} 
\end{equation}
where $\beta=0$ if $N-k\geq 1$, $\beta=1$ if $N-k=0$, and
 $\beta=\infty$ 
if $N-k<0$.
\end{theorem}
We can also compute $\tilde{L}^{N,k,d}_{m}$ only by using the above
 theorem,
 and this process is an analogue of the B-model computation in the 
Calabi-Yau case. Explicitly, the following recursive formula 
of the virtual structure constants holds.
\begin{cor}
The virtual structure constants $\tilde{L}_{n}^{N,k,d}$ can be fully 
reconstructed from the relation (\ref{hyper}).
As a result, we can compute all the virtual structure constants by
 using the   
initial condition and the recursive formula:
\begin{eqnarray}
&&\sum_{n=0}^{k-1}\tilde{L}_{n}^{N,k,1}w^{n}=
k\cdot\prod_{j=1}^{k-1}(jw+(k-j)), \no\\
&&\sum_{m=0}^{N-1+(k-N)d}\tilde{L}_{m}^{N,k,d}z^{m}=\no\\
&&\sum_{l=2}^{d}(-1)^{l}\sum_{0=i_{0}<\cdots<i_{l}=d}
\sum_{j_{l}=0}^{N-1+(k-N)d}\cdots\sum_{j_{2}=0}^{j_{3}}\sum_{j_{1}=0}^{j_{2}}
\prod_{n=1}^{l}\biggl((1+(d-i_{n-1})(\frac{z-1}{d}))^{j_{n}-j_{n-1}}
\cdot \tilde{L}_{j_{n}+(N-k)i_{n-1}}^{N,k,i_{n}-i_{n-1}}\biggr).\no\\
\label{cinic}
\end{eqnarray}
\end{cor}
 We can regard (\ref{cinic}) as an alternate 
definition of the virtual structure constants. 
\subsection{Residue Integral Representation of Virtual Structure
 Constants} 
By solving the recursive formula (\ref{cinic}) for low degrees
 explicitly, we reached a residue integral 
representation of the virtual structure constants. In the following, we
 give some definitions necessary to describe 
the formula we have obtained. First, we define rational functions
 $F_{d}(z,w)\;(d\in{\bf N})$ in $z,w$ by,
\begin{defi}
\begin{eqnarray}
F_{d}(z,w)&:=&k\prod_{j=1}^{d-1}\bigl(\frac{d}{jz+(d-j)w}\bigr)^{N}
 \prod_{j=1}^{kd-1}\bigl(\frac{jz+(kd-j)w}{d}\bigr).
\end{eqnarray}
\end{defi}
Next, we introduce the ordered partition of a positive integer $d$, which
 plays a central role in this paper. 
\begin{defi}
Let $OP_{d}$ be the set of ordered partitions of the positive integer $d$:
\begin{equation}
OP_{d}=\{\sigma_{d}=(d_{1},d_{2},\cdots,d_{l(\sigma_{d})})\;\;|\;\;
\sum_{j=1}^{l(\sigma_{d})}d_{j}=d\;\;,\;\;d_{j}\in{\bf N}\}.
\label{part} 
\end{equation}
From now on, we denote a ordered partition $\sigma_{d}$ by 
$(d_{1},d_{2},\cdots,d_{l(\sigma_{d})})$. In (\ref{part}), we denote 
the length of the ordered partition $\sigma_{d}$ by $l(\sigma_{d})$.
\end{defi}
With this set up, the residue integral representation mentioned above
 is given as follows:
\begin{conj}
\begin{eqnarray}
&&\sum_{j=0}^{N-1-(N-k)d}\frac{\tilde{L}_{n}^{N,k,d}}{d}z^jw^{N-1-(N-k)d-j}\no\\
&&=\sum_{\sigma_{d}\in OP_{d}}
\frac{1}{(2\pi\sqrt{-1})^{l(\sigma_{d})-1}\prod_{j=1}^{l(\sigma_{d})}d_{j}}\oint_{C_{0}}
 dz_{1}\cdots
\oint_{C_{0}}
 dz_{l(\sigma_{d})-1}\prod_{j=1}^{l(\sigma_{d})-1}\frac{z_{j}^{1-N}}
{\biggl(
 \frac{z_{j}-z_{j-1}}{d_{j}}+\frac{z_{j}-z_{j+1}}{d_{j+1}}\biggr)}\prod_{j=1}^{l(\sigma_{d})}F_{d_{j}}(z_{j-1},z_{j}),\no\\
\label{int}
\end{eqnarray}
where $z_{0}=z,\;z_{l(\sigma_{d})}=w$.
\end{conj}
In (\ref{int}), $\frac{1}{2\pi\sqrt{-1}}\oint_{C_{0}}dz_{j}$ represents the
 operation of taking the residue at $z_{j}=0$. 
The residue integral in (\ref{int}) depends heavily on the order of
 integration, and we have to take residues of 
$z_{j}$'s in descending (or ascending) order of subscript $j$. 
We have a proof of the above formula up to $d=3$ and checked
 numerically its validity up to $d=6$.
Let us look at the formula (\ref{int}) more closely in the $d=1,2$
 cases.
In the $d=1$ case, the ordered partition of $1$ is just $(1)$, and
 $(\ref{int})$ reduces to,
\begin{equation}
\sum_{j=0}^{k-1}\tilde{L}_{n}^{N,k,1}z^jw^{k-1-j}=F_{1}(z,w)=
 k\prod_{j=1}^{k-1}(jz+(k-j)w),\no\\
\end{equation}
which is nothing but the initial condition of the virtual structure
 constants in (\ref{cinic}).
In the $d=2$ case, the ordered partitions of $2$ are $(2)$ and $(1,1)$, and
 (\ref{int}) takes the following form: 
\begin{eqnarray}
\sum_{j=0}^{N-1-2(N-k)}\frac{\tilde{L}_{n}^{N,k,2}}{2}z^jw^{N-1-2(N-k)-j}=\frac{1}{2}F_{2}(z,w)+\frac{1}{2\pi\sqrt{-1}}
\oint_{C_{0}}\frac{u^{1-N}du}{u-z+u-w}F_{1}(z,u)F_{1}(u,w).
\end{eqnarray}
\section{Geometrical Derivation}
In this section, we derive (\ref{int}) as an integral of Chern classes on
 the moduli space of polynomial maps with two marked points, by using the 
localization computation.
\subsection{Moduli Space of degree $d$ Polynomial Map with two marked
 Points and its Fixed Point Sets}
Let ${\bf a}_{j},\;(j=0,1,\cdots,d)$ be vectors in ${\bf C}^{N}$ 
and let $\pi_{N}: {\bf C}^{N}\rightarrow CP^{N-1}$ be the projection map.
 
In this paper, we define a degree $d$ polynomial map $p$ from ${\bf
 C}^{2}$ to ${\bf C}^{N}$ as the map that consists of 
${\bf C}^{N}$vector-valued degree $d$ homogeneous polynomials in two
 coordinates $s,t$ of ${\bf C}^{2}$:    
\begin{eqnarray}
&&p:{\bf C}^{2}\rightarrow {\bf C}^{N}\no\\
&&p(s,t)={\bf a}_{0}s^{d}+{\bf a}_{1}s^{d-1}t+{\bf
 a}_{2}s^{d-2}t^{2}+\cdots+{\bf a}_{d}t^{d}.
\label{polyp}
\end{eqnarray}
The map space is described by ${\bf C}^{N(d+1)}=\{ ({\bf a}_{0},{\bf
 a}_{1},\cdots,{\bf a}_{d}) \}$.
We denote by $Mp_{0,2}(N,d)$ the space obtained from dividing $\{({\bf
 a}_{0},\cdots,{\bf a}_{d})\in
{\bf C}^{N(d+1)}|\;{\bf a}_{0}\neq {\bf 0},\;{\bf a}_{d}\neq {\bf 0}\}$
 by two ${\bf C}^{\times}$ actions induced from the following 
two ${\bf C}^{\times}$ actions on ${\bf C}^{2}$ via the map $p$ in
 (\ref{polyp}). 
\begin{equation}
(s,t)\rightarrow ( \mu s,\mu t),\;\;(s,t)\rightarrow (s,\nu t).
\label{torus1}
\end{equation}
By the above two torus actions, $Mp_{0,2}(N,d)$ can be regarded as a
 parameter space of degree $d$ birational maps from $CP^{1}$ to $CP^{N-1}$
 with
two marked points in $CP^{1}$: $0(=(1:0))$ and $\infty(=(0:1))$ .
 In particular, the second torus action corresponds to the automorphism group of
 $CP^{1}$ that keeps $0$ and $\infty$ invariant.
The condition ${\bf a}_{0}, {\bf a}_{d}\neq {\bf 0}$ ensures that the
 images of $0$ and $\infty$ are well-defined in $CP^{N-1}$. 
In the $d=1$ case, $Mp_{0,2}(N,1)$ is identified with $CP^{N-1}\times
 CP^{N-1}$ by the two torus actions (\ref{torus1}). But in 
the $d\geq 2$ cases, $Mp_{0,2}(N,d)$ has singularities, and we have to
 resolve them. We denote by $\widetilde{Mp}_{0,2}(N,d)$ the space 
obtained after resolution. This $\widetilde{Mp}_{0,2}(N,d)$ is the
 moduli space of degree $d$ polynomial maps with two marked points.
 Let us consider the resolution of $Mp_{0,2}(N,2)$ as an example.
 $Mp_{0,2}(N,2)$ is obtained 
by dividing the space,
\begin{equation}
 \{({\bf a}_{0},{\bf a}_{1},{\bf a}_{2})\;|\;{\bf a}_{j}\in{\bf
 C}^{N},\;{\bf a}_{0}\neq {\bf 0},\;{\bf a}_{2}
\neq {\bf 0}\;\}
\end{equation}
 by the two torus actions: 
\begin{eqnarray}
&&({\bf a}_{0},{\bf a}_{1},{\bf a}_{2})\rightarrow (\mu^{2}{\bf
 a}_{0},\mu^{2}{\bf a}_{1},\mu^{2}{\bf a}_{2}),\no\\
&&({\bf a}_{0},{\bf a}_{1},{\bf a}_{2})\rightarrow ({\bf a}_{0},\nu{\bf
 a}_{1},\nu^{2}{\bf a}_{2}).
\label{torus1.2}
\end{eqnarray}
Since ${\bf a}_{0}\neq{\bf 0}$, we can use the first torus action to
 reduce $({\bf a}_{0},{\bf a}_{1},{\bf a}_{2})$ to
$(\pi_{N}({\bf a}_{0}),{\bf a}_{1},{\bf a}_{2})$. By the second torus
 action, the locus $(\pi_{N}({\bf a}_{0}),{\bf 0},{\bf a}_{2})$
becomes singular and we have to blow it up. Then the exceptional locus
 can be
 identified with $(\pi_{N}({\bf a}_{0}),\pi_{N}({\bf a}_{1}),{\bf
 a}_{2})$ 
with ${\bf a}_{1}\neq {\bf 0}$. At this locus, the second torus action
 acts in the following way:
\begin{equation}
(\pi_{N}({\bf a}_{0}),\pi_{N}({\bf a}_{1}),{\bf a}_{2})\rightarrow
 (\pi_{N}({\bf a}_{0}),\pi_{N}({\bf a}_{1}),\nu^{2}{\bf a}_{2}).
\end{equation}    
Since ${\bf a}_{2}\neq {\bf 0}$, the exceptional locus can be identified
 with 
$\{(\pi_{N}({\bf a}_{0}),\pi_{N}({\bf a}_{1}),\pi_{N}({\bf
 a}_{2}))\}=(CP^{N-1})^{3}$. This result suggests that we have to 
consider a chain of two degree 1 polynomial maps:
\begin{equation}
({\bf a}_{0}s_{1}+{\bf a}_{1}t_{1})\cup({\bf a}_{1}s_{2}+{\bf
 a}_{2}t_{2}),
\label{chain1}
\end{equation}
in addition to the usual degree $2$ polynomial maps. In (\ref{chain1}), the
 two torus actions (\ref{torus1}) are extended to 
each $(s_{j},t_{j}),\;(j=1,2)$. In the general $d$ case, the resolution of
 $Mp_{0,2}(N,d)$ forces us to consider a
chain of polynomial maps labeled by ordered partition
 $\sigma_{d}=(d_{1},d_{2},\cdots,d_{l(\sigma_{d})})$:
\begin{equation}
\cup_{j=1}^{l(\sigma_{d})}\bigl(\sum_{m_{j}=0}^{d_{j}} {\bf
 a}_{\sum_{i=1}^{j-1}d_{i}+m_{j}}(s_{j})^{m_{j}}(t_{j})^{d_{j}-m_{j}}\bigr),
\;\;\bigl({\bf a}_{\sum_{i=1}^{j}d_{i}}\neq {\bf
 0},\;\;j=0,1,\cdots,l(\sigma_{d})\bigr),
\label{chain2}
\end{equation}
where the two torus actions (\ref{torus1}) are extended to each
 $(s_{j},t_{j}),\;(j=1,2,\cdots,l(\sigma_{d}))$. 

From now on, we introduce the following ${\bf C}^{\times}$ action on
 ${\bf C}^{N(d+1)}$ and determine the fixed point set 
of $Mp_{0,2}(N,d)$.
\begin{equation}
({\bf a}_{0},{\bf a}_{1},\cdots,{\bf a}_{d}) \rightarrow
 (e^{\lambda_{0}t}{\bf a}_{0},e^{\lambda_{1}t}{\bf a}_{1},\cdots,e^{\lambda_{d}t}{\bf
 a}_{d}) 
\label{torus2}
\end{equation}
First, we look at the simplest map in $\widetilde{Mp}_{0,2}(N,d)$ given
 by,\footnote{If ${\bf a}_{0}={\bf a}_{d}$, 
the map (\ref{type1}) becomes constant map from $CP^{1}$ to $CP^{N-1}$,
 but we don't eliminate this locus in this paper.} 
\begin{equation}
{\bf a}_{0}s^{d}+{\bf a}_{d}t^{d}.
\label{type1}
\end{equation}
We can easily see that this map is invariant under the torus action
 (\ref{torus2}) because of the second torus action of 
(\ref{torus1}), hence we call the map (\ref{type1}) the type I fixed
 point.  
Due to the two torus actions in (\ref{torus1}), the type I fixed point set
 is identified with
 $ ( \pi_{N}({\bf a}_{0}) ,\pi_{N}({\bf a}_{d})) \in (CP^{N-1})^{2}$.
We also have to note that ${\bf Z}_{d}=\{ \exp(\frac{2\pi\sqrt{-1}
 j}{d}) |\;j=0,1,\cdots, d-1\}$ 
naturally acts on $CP^{1}$ in the following way:
\begin{equation}
(s:t)\rightarrow (s:\exp(\frac{2\pi\sqrt{-1} j}{d})t) 
\end{equation}
 but leaves $( \pi_{N}({\bf a}_{0}) ,\pi_{N}({\bf a}_{d}))$ invariant.
 This means that when integrating over this fixed point set, we have to
 divide
the result by $|{\bf Z}_{d}|=d$. 

Let us consider the normal bundle of this fixed point set in
 $\widetilde{Mp}_{0,2}(N,d)$. Obviously, this bundle is spanned by the degrees of freedom
 of 
deformation of the map (\ref{type1}) by using ${\bf a}_{j},\;(j=1,2,\cdots,
 d-1)$. Therefore, normal vector space is given by the following 
$N(d-1)$ dimensional vector space:
\begin{equation}
{\bf C}^{N}s^{d-1}t\oplus{\bf C}^{N}s^{d-2}t^2\oplus{\bf
 C}^{N}s^{d-3}t^3\oplus\cdots\oplus{\bf C}^{N}st^{d-1}.
\label{normal1}
\end{equation}
We can easily see that the fixed points coming from usual degree $d$
 polynomial maps are exhausted by type I
fixed points. The remaining fixed points can be found from the exceptional
 locus given by the chain of polynomial maps (\ref{chain2}).
In (\ref{chain2}), we have extended the two torus actions (\ref{torus1})
 to each $(s_{j},t_{j})$.   
By using this fact, we can construct a chain of type I-like maps which
 remains fixed under the torus action (\ref{torus2}):   
\begin{eqnarray}
\cup_{j=1}^{l(\sigma_{d})} ({\bf
 a}_{\sum_{i=1}^{j-1}d_{i}}(s_{j})^{d_{j}}+{\bf a}_{\sum_{i=1}^{j}d_{i}}(t_{j})^{d_{j}}).
\label{fix2}
\end{eqnarray}
Note that image of the chain of maps (\ref{fix2}) in $CP^{N-1}$ is a nodal
 rational curve with $l(\sigma_{d})-1$ nodal singularities given by 
$\pi_{N}({\bf
 a}_{\sum_{i=1}^{j}d_{i}}),\;(j=1,2,\cdots,l(\sigma_{d})-1)$.  

We call the chain of maps given in (\ref{fix2}) a type II fixed point.
The type II fixed point set labeled by $\sigma_{d}$ is identified with 
$$\bigl(\pi_{N}({\bf a}_{0}),\pi_{N}({\bf a}_{d_{1}}),\pi_{N}({\bf
 a}_{d_{1}+d_{2}}), \cdots,
\pi_{N}({\bf a}_{d-d_{l(\sigma_{d})}}) ,\pi_{N}({\bf a}_{d})
 \bigr)\in(CP^{N-1})^{l(\sigma_{d})+1}.$$
Let $p_{j}:CP^{1}\rightarrow CP^{N-1}$ be $\pi_{N}({\bf
 a}_{\sum_{i=1}^{j-1}d_{i}}(s_{j})^{d_{j}}+{\bf
 a}_{\sum_{i=1}^{j}d_{i}}(t_{j})^{d_{j}})$ and $C_{j}$ be the image of $p_{j}$ in 
$CP^{N-1}$.
In the same way as in the type I case, ${\bf Z}_{d_{j}}$  acting on
 $(s_{j},t_{j})$ keeps $p_{j}$ invariant. Therefore, 
in integrating over the type II fixed point set labeled by $\sigma_{d}$, we
 have to divide
the result by $\prod_{j=1}^{l(\sigma_{d})}|{\bf
 Z}_{d_{j}}|=\prod_{j=1}^{l(\sigma_{d})}d_{j}$.  

The normal vector space of the type II fixed point set labeled by $\sigma_{d}$
 is spanned by the degrees of freedom coming from deforming $p_{j}$ individually and 
by the degrees of freedom associated to resolving nodal singularities of the image
 curve: 
\begin{eqnarray}
&&\oplus_{j=1}^{l(\sigma_{d})}\biggl({\bf
 C}^{N}s_{j}^{d_{j}-1}t_{j}\oplus{\bf C}^{N}s_{j}^{d_{j}-2}t_{j}^2\oplus{\bf
 C}^{N}s_{j}^{d_{j}-3}t_{j}^3
\oplus\cdots\oplus{\bf C}^{N}s_{j}t_{j}^{d_{j}-1}\biggr)\no\\
&&\oplus_{j=1}^{l(\sigma_{d})-1}\biggl(T^{\prime}_{\infty}C_{j}\otimes
 T^{\prime}_{\bf 0}C_{j+1}\biggr).
\label{norm2}
\end{eqnarray}
\subsection{Localization Computation of Virtual Structure Constants}
In this subsection, we compute the $\frac{\tilde{L}^{N,k,d}_{n}}{d}$, which
 may be regarded as the B-model analogue of 2-pointed Gromov-Witten
 invariants: 
\begin{eqnarray}
&&\frac{1}{dk}\langle{\cal O}_{h^{N-2-n}}{\cal O}_{h}{\cal
 O}_{h^{n-1+(N-k)d}}\rangle_{0,d}
 =\frac{1}{k}\langle{\cal O}_{h^{N-2-n}}{\cal
 O}_{h^{n-1+(N-k)d}}\rangle_{0,d}\no\\
&&=\frac{1}{k}\int_{\overline{M}_{0,2}(CP^{N-1},d)}c_{top}\bigl(R^{0}(\pi_{*}ev_{3}^{*}({\cal
 O}(k)))\bigr)\wedge ev_{1}^{*}(h^{N-2-n})
\wedge ev_{2}^{*}(h^{n-1+(N-k)d}),
\label{corr}
\end{eqnarray}
by using the moduli space introduced in the previous subsection. 
In (\ref{corr}), $h$ is the hyperplane class of $CP^{N-1}$, $
 \overline{M}_{0,n}(CP^{N-1},d) $ represents moduli space of degree $d$ 
stable maps from genus $0$ stable curve to $CP^{N-1}$ with $n$ marked
 points,  
$ev_{i}: \overline{M}_{0,n}(CP^{N-1},d) \rightarrow CP^{N-1}$ is the
 evaluation map of  the $i$-th marked point and
 $\pi:\overline{M}_{0,3}(CP^{N-1},d) 
\rightarrow \overline{M}_{0,2}(CP^{N-1},d)$ is the forgetful map.  
We also use the localization technique (Bott residue formula) 
associated with the torus action (\ref{torus2}). However, in the following
 computation, we only consider the case $\lambda_{0}=\lambda_{1}=
\cdots=\lambda_{d}=0$, for simplicity. To compensate for this choice,
 we have to treat the order of integration carefully. We will discuss 
these subtleties in the last part of this subsection.

First, we determine the contribution from the type I fixed point set
 $(CP^{N-1})^{2}$. From now on, we denote by $(CP^{N-1})_{0}$ (resp. 
$(CP^{N-1})_{1}$) the first (resp. the second) $CP^{N-1}$ of
 $(CP^{N-1})^{2}$.   
We define $h_{i}$ as the hyperplane class of $(CP^{N-1})_{i}$.  Let
 $p_{1}:CP^{1}\rightarrow CP^{N-1}$ be the map defined by,
\begin{equation}
p_{1}(s:t):=\pi_{N}({\bf a}_{0}s^{d}+{\bf a}_{d}t^{d}).
\label{map1}
\end{equation} 
In the construction of $\widetilde{Mp}_{0,2}(N,d)$, the two marked points
 of $CP^{1}$ are fixed to $0=(1:0)$ and $\infty=(0:1)$. Obviously, 
we have 
\begin{equation}
p_{1}(1:0)=\pi_{N}({\bf a}_{0}),\;\;p_{1}(0:1)=\pi_{N}({\bf a}_{d}).
\end{equation}
Therefore, the classes that correspond to $ev_{1}^{*}(h^{N-2-n})$ and
 $ev_{2}^{*}(h^{n-1+(N-k)d})$ are given by 
$h_{0}^{N-2-n}$ and $h_{1}^{n-1+(N-k)d}$ respectively. Then we consider
 the vector bundle corresponding to 
$R^{0}\pi_{*}ev_{3}^{*}{\cal O}(k)$. The corresponding vector space is
 spanned by $H^{0}(CP^{1}, p_{1}^{*}{\cal O}(k))$, and can written as 
\begin{equation}
\oplus_{j=0}^{kd}{\bf C}s^{j}t^{kd-j}.
\end{equation}
In the localization computation, we can identify $s$ with ${\cal
 O}_{(CP^{N-1})_{0}}(\frac{1}{d})$ and $t$ with 
${\cal O}_{(CP^{N-1})_{1}}(\frac{1}{d})$ through (\ref{map1}).
 Therefore, the first Chern class of ${\bf C}s^{i}t^{j}$ is given by,
\begin{equation}
\frac{ih_{0}+jh_{1}}{d}.
\end{equation}
In this way, the class that corresponds to
 $c_{top}(R^{0}\pi_{*}ev_{3}^{*}{\cal O}(k))$ turns out to be 
\begin{equation}
\prod_{j=0}^{kd}\biggl(\frac{jh_{0}+(kd-j)h_{1}}{d}\biggr).
\end{equation} 
If we look back at (\ref{normal1}), we can also determine the top Chern
 class of the normal bundle of the type I fixed point set as follows:
\begin{equation}
\prod_{j=1}^{d-1}\biggl(\frac{jh_{0}+(d-j)h_{1}}{d}\biggr)^{N}.
\end{equation}
Putting these pieces together, we can write down the contribution from
the type I fixed point set by the localization theorem: 
\begin{equation}
\frac{1}{d}\int_{(CP^{N-1})^2}h_{0}^{N-2-n}h_{1}^{n-1+(N-k)d}
 \frac{\prod_{j=0}^{kd}(\frac{jh_{0}+(kd-j)h_{1}}{d})}{\prod_{j=1}^{d-1}
(\frac{jh_{0}+(d-j)h_{1}}{d})^{N}},
\label{loc1}
\end{equation}
where the factor$\frac{1}{d}$ comes from the ${\bf Z}_{d}$ action
 mentioned in the previous subsection.

Next, we determine the contribution from the type II fixed point set
 labeled by $\sigma_{d}$. 
 Let $(CP^{N-1})_{i}$ be the $(i+1)$-th $CP^{N-1}$ of
 $(CP^{N-1})^{l(\sigma_{d})+1}$ considered as a type II fixed point set, and 
let $h_{i}$ be its hyperplane class.
The computation goes in the same way as in the type I case, except for
the effect of nodal singularities. Therefore, we consider here 
the contributions of these singularities. As for the
 vector bundle corresponding to 
$R^{0}\pi_{*}ev_{3}^{*}{\cal O}(k)$, we have to consider the exact
 sequence:
\begin{eqnarray}
0\rightarrow
 H^{0}(\cup_{j=1}^{l(\sigma_{d})}C_{j},(\cup_{j=1}^{l(\sigma_{d})}p_{j})^{*}{\cal O}(k))
\rightarrow \oplus_{j=1}^{l(\sigma_{d})}H^{0}(C_{j},p_{j}^{*}{\cal
 O}(k))\rightarrow
\oplus_{j=1}^{l(\sigma_{d})-1}{\cal O}_{\pi_{N}({\bf
 a}_{\sum_{i=1}^{j}d_{i}})}(k)\rightarrow 0.
\label{exact1}
\end{eqnarray} 
From this exact sequence, we can easily see that we have to insert an additional
 $\prod_{j=1}^{l(\sigma_{d})-1}\frac{1}{kh_{j}}$.
We next turn to the effect of $T'_{\infty}C_{j}\otimes T'_{\bf
 0}C_{j+1}$ in (\ref{norm2}).
Obviously this can be written as,
\begin{equation}
{\bf C} \frac{d}{d(\frac{s_{j}}{t_{j}})} \otimes
  \frac{d}{d(\frac{t_{j+1}}{s_{j+1}})},
\end{equation}
and the $s_{j}$, $t_{j}$ are identified with ${\cal
 O}_{(CP^{N-1})_{j-1}}(\frac{1}{d_{j}})$, ${\cal O}_{(CP^{N-1})_{j}}(\frac{1}{d_{j}})$
respectively. Hence its first Chern class is given by,
\begin{equation}
\frac{h_{j}-h_{j-1}}{d_{j}}+\frac{h_{j}-h_{j+1}}{d_{j+1}}.
\end{equation} 
By the localization theorem, we also have to insert 
$\prod_{j=1}^{l(\sigma_{d})-1}\bigl(\frac{h_{j}-h_{j-1}}{d_{j}}+\frac{h_{j}-h_{j+1}}{d_{j+1}}\bigr)^{-1}$.
Combining these considerations, we can write down the contributions from
 the type II fixed point set labeled by $\sigma_{d}$:
\begin{eqnarray}
&& \frac{1}{\prod_{j=1}^{l(\sigma_{d})}d_{j}}
 \int_{(CP^{N-1})^{l(\sigma_{d})+1}}h_{0}^{N-2-n}h_{l(\sigma_{d})}^{n-1+(N-k)d}
\biggl(\prod_{j=1}^{l(\sigma_{d})}\frac{\prod_{i=0}^{kd_{j}}(\frac{ih_{j-1}+(kd_{j}-i)h_{j}}{d_{j}})}{\prod_{i=1}^{d_{j}-1}
(\frac{ih_{j-1}+(d_{j}-i)h_{j}}{d_{j}})^{N}}
\biggr)\prod_{j=1}^{l(\sigma_{d})-1}
 \frac{1}{kh_{j}\bigl(\frac{h_{j}-h_{j-1}}{d_{j}}+\frac{h_{j}-h_{j+1}}{d_{j+1}}\bigr) },\no\\
\label{loc2} 
\end{eqnarray}
where the factor $\frac{1}{\prod_{j=1}^{l(\sigma_{d})}d_{j}}$ comes
 from the ${\bf Z}_{d_{j}}$ action on $C_{j}$
 $(j=1,2,\cdots,l(\sigma_{d}))$. 
Then we integrate out (\ref{loc1}) and (\ref{loc2}) and divide the
 result by $k$.  After all this, 
we arrive at the corresponding summands in the formula (\ref{int}).
 Lastly, we have to mention order of integration. In practice, 
we have to order the integrations of all the summands in
 descending (or ascending) order of subscript $j$ of $(CP^{N-1})_{j}\subset
(CP^{N-1})^{l(\sigma_{d})+1}$.  
\section{Applications to Local Mirror Symmetry of  Vector Bundles over
 $CP^{1}$}
Our geometrical computation is also applicable to complete intersections
 in $CP^{N-1}$ and to the local geometries 
$\oplus_{i=1}^{m}{\cal O}(k_{i})\rightarrow CP^{N-1}$. In this section,
 we consider two examples of local mirror symmetry of 
$CP^{1}$, ${\cal O}(-1)\oplus{\cal O}(-1)\rightarrow CP^{1}$ and ${\cal
 O}(1)\oplus{\cal O}(-3)\rightarrow CP^{1}$.
The first one is the simplest example of local mirror symmetry, and the
 second one is a typical example of non-nef local mirror 
symmetry, which was analyzed extensively in \cite{fj1},\cite{fj2}.
\subsection{${\cal O}(-1)\oplus{\cal O}(-1)\rightarrow CP^{1}$}
In this example, we compute the virtual structure constants $\alpha_{n}^{d}$
 that correspond to local Gromov-Witten invariants:
\begin{equation}
\langle{\cal O}_{h^{n}}{\cal
 O}_{h^{2-n}}\rangle_{0,d}=\int_{\overline{M}_{0,2}(CP^{1},d)}
\bigl(c_{top}(R^{1}\pi_{*}ev_{3}^{*}{\cal O}(-1))\bigr)^{2}\wedge
 ev_{1}^{*}(h^{n})\wedge ev_{2}^{*}(h^{2-n}),\;\;(n=0,1,2).
\end{equation} 
Following the computation in the previous section, we obtain a closed
 formula that computes $\alpha_{n}^{d}$:
\begin{eqnarray}
\alpha_{n}^{d}=\sum_{\sigma_{d}\in
 OP_{d}}\frac{1}{\prod_{j=1}^{l(\sigma_{d})}d_{j}}\int_{(CP^{1})^{l(\sigma_{d})+1}}
h_{0}^{n}h_{l(\sigma_{d})}^{2-n}\prod_{j=1}^{l(\sigma_{d})}G_{d_{j}}(h_{j-1},h_{j})\prod_{j=1}^{l(\sigma_{d})-1}
\frac{(-h_{j})^{2}}{\frac{h_{j}-h_{j-1}}{d_{j}}+\frac{h_{j}-h_{j+1}}{d_{j+1}}},
\label{-1-1}
\end{eqnarray}
where
\begin{equation}
G_{d}(x,y):=\frac{\prod_{j=1}^{d-1}\bigl(\frac{-jx-(d-j)y}{d}\bigr)^{2}}
{\prod_{j=1}^{d-1}\bigl(\frac{jx+(d-j)y}{d}\bigr)^{2}}=1.
\label{simple}
\end{equation}
In deriving (\ref{-1-1}), the exact sequence (\ref{exact1}) is replaced
 by the following exact sequence:
\begin{eqnarray}
0\rightarrow \oplus_{j=1}^{l(\sigma_{d})-1}{\cal O}_{\pi_{N}({\bf
 a}_{\sum_{i=1}^{j}d_{i}})}(-1)\rightarrow 
H^{1}(\cup_{j=1}^{l(\sigma_{d})}C_{j},(\cup_{j=1}^{l(\sigma_{d})}p_{j})^{*}{\cal
 O}(-1))
\rightarrow \oplus_{j=1}^{l(\sigma_{d})}H^{1}(C_{j},p_{j}^{*}{\cal
 O}(-1))\rightarrow 0.
\label{exact2}
\end{eqnarray} 
Since (\ref{simple}) holds,  (\ref{-1-1}) can be further simplified to, 
\begin{eqnarray}
\alpha_{n}^{d}=\sum_{\sigma_{d}\in
 OP_{d}}\frac{1}{\prod_{j=1}^{l(\sigma_{d})}d_{j}}\int_{(CP^{1})^{l(\sigma_{d})+1}}
h_{0}^{n}h_{l(\sigma_{d})}^{2-n}\prod_{j=1}^{l(\sigma_{d})-1}
\frac{(h_{j})^{2}}{\frac{h_{j}-h_{j-1}}{d_{j}}+\frac{h_{j}-h_{j+1}}{d_{j+1}}}.
\label{-1-1s}
\end{eqnarray}
Due to the fact that $h_{j}^{2}=0$ in $H^{*}((CP^{1})_{j},{\bf C})$,
 the summand in (\ref{-1-1s}) vanishes if $l(\sigma_{d})>1$. Hence we
 obtain
\begin{equation}
\alpha_{n}^{d}=\frac{1}{d}\delta_{n1}.
\label{re-1}
\end{equation}
If we consider this in analogy with the Calabi-Yau hypersurface case
 \cite{cj}, we expect that 
$t=x+\sum_{d=1}^{\infty}\alpha_{0}^{d}e^{dx}$ gives us the mirror map
 of this model. (\ref{re-1}) says that it
 is trivial in this case. Therefore, $\alpha_{1}^{d}$ 
should coincide with the corresponding local Gromov-Witten invariant
 $\langle{\cal O}_{h}{\cal O}_{h}\rangle_{0,d}$.
(\ref{re-1}) agrees with this expectation.
\subsection{${\cal O}(1)\oplus{\cal O}(-3)\rightarrow CP^{1}$}
In this case, we compute equivariant virtual structure constants
 $\beta_{mn}^{d}(z)$ that correspond to the 
equivariant Gromov-Witten invariants:
\begin{equation}
\langle{\cal O}_{h^{m}}{\cal
 O}_{h^{n}}(z)\rangle_{0,d}=\int_{\overline{M}_{0,2}(CP^{1},d)}
\biggl(\frac{ \sum_{j=0}^{3d-1}c_{j}(R^{1}\pi_{*}ev_{3}^{*}{\cal
 O}(-3)) }
{\sum_{j=0}^{d}z^{j}c_{j}(R^{0}\pi_{*}ev_{3}^{*}{\cal
 O}(1))}\biggr)\wedge ev_{1}^{*}(h^{m})\wedge ev_{2}^{*}(h^{n}),\;(0\leq m,n\leq 1).
\end{equation}
In the same way as in the previous subsection, $\beta_{mn}^{d}(z)$ is
 given as follows:
\begin{eqnarray}
\beta_{mn}^{d}(z)=\sum_{\sigma_{d}\in
 OP_{d}}\frac{1}{\prod_{j=1}^{l(\sigma_{d})}d_{j}}\int_{(CP^{1})^{l(\sigma_{d})+1}}
h_{0}^{m}h_{l(\sigma_{d})}^{n}\prod_{j=1}^{l(\sigma_{d})}H_{d_{j}}(h_{j-1},h_{j},z)\prod_{j=1}^{l(\sigma_{d})-1}
\frac{(1+zh_{j})(1-3h_{j})}{\frac{h_{j}-h_{j-1}}{d_{j}}+\frac{h_{j}-h_{j+1}}{d_{j+1}}},
\label{-3s}
\end{eqnarray}
where
\begin{equation}
H_{d}(x,y,z):=\frac{\prod_{j=1}^{3d-1}\bigl(1-\frac{jx+(3d-j)y}{d}\bigr)}
{\prod_{j=0}^{d}(1+z\frac{jx+(d-j)y}{d})\prod_{j=1}^{d-1}\bigl(\frac{jx+(d-j)y}{d}\bigr)^{2}}.
\end{equation}
Let us compute $\beta_{00}^{d}(z),\beta_{10}^{d}(z),\beta_{11}^{d}(z)$
 by using the formula (\ref{-3s}). Here, we show the results 
for lower degrees by using the generating function
 $\beta_{ij}(e^{x},z):=\sum_{d=1}^{\infty}\beta_{ij}^{d}(z)e^{dx}$.
\begin{eqnarray}
\beta_{00}(e^{x},z)&:=&(6z + z^{2} + 5)e^{x} + ({\displaystyle \frac {
645}{4}} z + {\displaystyle \frac {311}{4}} z^{2} + 
{\displaystyle \frac {63}{4}} z^{3} + {\displaystyle \frac {5}{
4}} z^{4} + 104)e^{2x} \no\\
&&+ (6387z + {\displaystyle \frac {91421}{18}} z^{2} + 
{\displaystyle \frac {8767}{18}} z^{4} + {\displaystyle \frac {
4197}{2}} z^{3} + {\displaystyle \frac {83723}{27}}  + 
{\displaystyle \frac {121}{2}} z^{5} + {\displaystyle \frac {85
}{27}} z^{6})e^{3x}+\cdots\no\\ 
\beta_{10}(e^{x},z)&:=&( - 3 - z)e^{x} + ( - {\displaystyle \frac
 {139}{4}
} z - {\displaystyle \frac {33}{4}} z^{2} - {\displaystyle 
\frac {3}{4}} z^{3} - {\displaystyle \frac {177}{4}} )e^{2x}
 \no\\
&& + ( - 1131 - {\displaystyle \frac {5917}{4}} z - 762z
^{2} - 28z^{4} - {\displaystyle \frac {407}{2}} z^{3} - 
{\displaystyle \frac {19}{12}} z^{5})e^{3x}+\cdots\no\\
\beta_{11}(e^{x},z)&:=&e^{x} + ({\displaystyle \frac {31}{2}}  + 
{\displaystyle \frac {9}{2}} z + {\displaystyle \frac {1}{2}} 
z^{2})e^{2x} + (380 + {\displaystyle \frac {549}{2}} z + 
{\displaystyle \frac {175}{2}} z^{2} + {\displaystyle \frac {5
}{6}} z^{4} + {\displaystyle \frac {27}{2}} z^{3})e^{3x}+\cdots
\end{eqnarray}
From these results, we can see that $\beta_{10}(e^{x},z)$ coincides
 with 
the equivariant mirror map $\tilde{t}(q,\lambda)$ in the formula (3.23)
 of \cite{fj1}. 
Moreover, if we compute $\beta_{00}(e^{x},z)+(z-3)\beta_{10}(e^{x},z)$\footnote{This combination is derived 
by introducing the metric $\eta_{ij}$ induced from classical intersection number $\eta_{ij}:=\int_{CP^{1}}\frac{h^{i+j}}{(1+zh)(1-3h)}$.},
 the result 
turns out to be, 
\begin{eqnarray}
&&(14 + 6z)e^{x} + ({\displaystyle \frac {39}{4}} z^{3} + 
{\displaystyle \frac {271}{4}} z^{2} + {\displaystyle \frac {1
}{2}} z^{4} + {\displaystyle \frac {885}{4}} z + 
{\displaystyle \frac {947}{4}} )e^{2x} \no\\
&&+({\displaystyle \frac {38775}{4}} z + {\displaystyle \frac {
211885}{36}} z^{2} + {\displaystyle \frac {149}{4}} z^{5} + 
{\displaystyle \frac {3308}{9}} z^{4} + {\displaystyle \frac {
169}{108}} z^{6} + {\displaystyle \frac {175334}{27}}  + 1947
z^{3})e^{3x}+\cdots, 
\end{eqnarray}
which is nothing but the other equivariant mirror map $t(q,\lambda)$ in
 the formula (3.23) of \cite{fj1}! 
In this way, we have computed full equivariant mirror map  without 
using Birkhoff factorization inevitable in the analysis in \cite{fj1}.
 Then we set 
\begin{eqnarray}
&&\tilde{L}_{1}(e^{x},z):=1+\partial_{x}(\beta_{00}(e^{x},z)+(z-3)\beta_{10}(e^{x},z))=\frac{\partial
 t}{\partial{x}},\no\\
&& \tilde{L}_{2}(e^{x},z) :=\partial_{x}\beta_{11}(e^{x},z),
\label{-3m}
\end{eqnarray}
in analogy with the Calabi-Yau hypersuface case \cite{cj}. The
 first line of (\ref{-3m}) also asserts,
\begin{eqnarray}
t(x,z)=x+\beta_{00}(e^{x},z)+(z-3)\beta_{10}(e^{x},z).
\label{mirror}
\end{eqnarray}
If we invert (\ref{mirror}) regarding $z$ as a parameter, we obtain
 $x(t,z)$. Then again by using the same analogy as before,
we compute
 $\frac{\tilde{L}_{2}(e^{x(t,z)},z)}{\tilde{L}_{1}(e^{x(t,z)},z)}$. The result turns out to be,
\begin{equation}
e^{t} + ( - 3z + 3 + z^{2})e^{2t} + (
{\displaystyle \frac {69}{4}} z^{2} - {\displaystyle \frac {81
}{4}} z + {\displaystyle \frac {39}{4}}  + z^{4} - 
{\displaystyle \frac {27}{4}} z^{3})e^{3t}+\cdots, 
\end{equation}
which is nothing but the equivariant A-model Yukawa coupling compatible
 with the formula (3.23) of \cite{fj1}!!
We can also apply this technique to the model ${\cal O}(1) \oplus {\cal
 O}(-1) \oplus {\cal O}(-1) \oplus {\cal O}(-1)
\rightarrow CP^{1}$ and obtain the same mirror map and A-model Yukawa
 coupling as the ones computed in \cite{fj2}. 
\section{Conclusion}
In this paper, we derived the virtual structure constants by applying
the localization computation 
to the moduli space of polynomial maps with two marked points. The process
 of computation is very 
similar to the well-known result of Kontsevich \cite{kont}, but is
 much simpler because we consider a 
simple moduli space of polynomial maps instead of the moduli space of
 stable maps. As a result, we 
obtain a B-model analogue of standard Gromov-Witten invariants. Unlike
 the standard Gromov-Witten 
invariants, the virtual structure constants need not vanish when they
 have an insertion of an operator
induced from the identity element of the cohomology ring. Moreover, the virtual
 structure constants with insertions of the 
identity operator give us the expansion coefficients of the mirror map in the
 examples we have treated. 
Therefore, our computation provides a geometrical construction of the
 mirror map. 

The line of thought in this paper 
stems from our endeavor to interpret geometrically the computation process
 of the generalized mirror transformation
\cite{gene}. What we are aiming at is to describe the generalized mirror
 transformation as a process of 
changing the moduli space of Gauged Linear Sigma Model into the one of
 stable maps (see the discussion in 
\cite{bert} in the nef cases) \cite{geo}. To this end, we need to
 characterize the expansion coefficients of 
mirror map geometrically, and we think that the construction given in
 this paper provides what we 
need. We also hope to generalize our 
construction to various weighted projective spaces, especially in the
 case when they have several K\"ahler 
forms. If we accomplish this task, we will obtain a concrete geometrical
 understanding of the mirror computation.  

\end{document}